\documentclass[10pt]{article}

\usepackage{amsmath,latexsym,amsthm,amsfonts,amssymb}

\newtheorem{theorem}{Theorem}
\newtheorem{problem}{Problem}
\newtheorem{observation}{Observation}
\newtheorem{proposition}{Proposition}
\newtheorem{corollary}{Corollary}
\newtheorem{definition}{Definition}
\newtheorem{lemma}{Lemma}
\begin{document}
 
\title{A semifilter approach to selection principles}
\author{Lubomyr Zdomsky}
\maketitle
\baselineskip15pt
\begin{abstract}
In this paper we develop the semifilter approach to the classical Menger and Hurewicz
properties and show that the small cardinal $\mathfrak g $ is a
 lower bound of the additivity number of the $\sigma$-ideal generated by
 Menger subspaces of the Baire space, and under $\mathfrak u<\mathfrak g$ 
every subset $X$ of the real line with the property $\mathrm{Split}(\Lambda,\Lambda)$ is Hurewicz, and thus it is consistent
with ZFC that the property $\mathrm{Split}(\Lambda,\Lambda)$ is preserved by unions of less than $\mathfrak b$
subsets of the real line.
\end{abstract}

\large \centerline{\textbf{Introduction}} \normalsize 
\smallskip

In this paper we shall present two directions of applications of semifilters
in selection principles on topological spaces. First, we shall consider preservation
by unions of the Menger property. 

Trying to describe the $\sigma$-compactness in terms of open covers, K.Menger introduced
in \cite{Me1} the following property, called \emph{ the Menger property}: a topological space
$X$ is said to have this  property if for every sequence $(u_n)_{n\in\omega}$
of open covers of $X$ there exists a sequence $(v_n)_{n\in\omega}$ such that
each $v_n$ is a finite subfamily of $u_n$ and the collection $\{\cup v_n:n\in\omega\}$
is a cover of $X$. The class of Menger topological spaces, i.e. spaces having the Menger
property appeared to be much wider than the class of $\sigma$-compact spaces (see \cite{BT}, \cite{CP}, \cite{JMSS} and many others),
but it has  interesting properties itself and poses a number of open questions.
One of them, namely the question about the value of additivity of corresponding $\sigma$-ideal,
will be discussed in this paper. Let us recall that a collection $\mathcal I$
of subsets of a set $X$ is called a \emph{$\sigma$-ideal} if it is closed under
taking subsets and countable unions. Therefore, the union $\cup\mathcal J$
belongs to $\mathcal I$ for every countable subfamily $\mathcal J$ of $\mathcal I$.
In light of this property of $\sigma$-ideals it is interesting to find the smallest
cardinality $\tau$ such that the union $\cup\mathcal J$ in not in $\mathcal I$
for some $\mathcal J\subset\mathcal I$ with $|\mathcal J|=\tau$.
When $\cup\mathcal I=X$ and $X\not\in\mathcal I$ such the cardinality 
obviously exists and we denote it by  $\mathrm{add}(\mathcal I)$.
 It is easy to prove (see, for example, \cite{JMSS})
that the collection $M(X)$ of subspaces of a topological space $X$ contained in subspaces
with the Menger property form a $\sigma$-ideal, so one can ask about the value of $\mathrm{add}(M(X))$. 
According to \cite{BST}, for the Baire  space $\mathbb{N}^\omega$ this additivity number
is situated between cardinals $\mathfrak b$ and $\mathrm{cf}(\mathfrak d)$,
where $\mathfrak b$ and $\mathfrak d$ are well-known bounding and dominating numbers 
respectively, see \cite{Va}. It was also conjectured in \cite{BST}
that $\mathrm{add}(M(X))=\mathfrak b$, see Problem 2.4 there.
We shall prove here that another small cardinal, namely $\mathfrak g$,
is an lower bound of $\mathrm{add}(M(X))$ for each hereditarily Lindel\"of
topological space $X$. Since there are models of ZFC with $\mathfrak b<\mathfrak g$
(see \cite{Va}),
this answers the above mentioned problem in negative.
 Concerning topological spaces which contain
non-Lindel\"of subspaces, the straightforward proof of the fact that this additivity  equals $\aleph_1$
is left to the reader.  

Another direction is devoted to splittability of open covers. Following \cite {GN} and \cite{JMSS}
we say that a family  $u$ of subsets of a  set $X$ is
\begin{itemize}
\item \emph{a large cover} of $X$, if every $x\in X$ belongs to infinitely many $U\in u$;
\item \emph{an $\omega$-cover}, if for every finite subset $K$ of $X$ the family
$\{U\in u:K\subset U\}$ is infinite;
\item \emph{a $\gamma$-cover}, if for every $x\in X$ the family $\{U\in u:x\not\in U\}$ is finite.
\end{itemize} 
From now on we denote by $\Lambda(X)$ (resp. $\Omega(X)$, $\Gamma(X)$) the family 
of all large (resp. $\omega$-, $\gamma$-) covers of $X$. A topological space $X$
satisfies the selection hypothesis $\mathrm{Split}(\Lambda,\Lambda)$, if for every $u\in\Lambda(X)$ there
are $v_1,v_2\in\Lambda(X)$ such that $v_1\cap v_2=\emptyset$ and $v_1\cup v_2\subset u$. 
The class $\mathrm{Split}(\Lambda,\Lambda)$ contains all Hurewicz spaces and all spaces with the Rothberger property,
see \cite[Cor. 29, Th. 15]{Sch}. Recall, that a topolgical space $X$ has the
 Hurewicz property, if for every sequence $(u_n)_{n\in\omega}$
of open covers of $X$ there exists a $\gamma$-cover $\{B_n:n\in\omega\}$ of
 $X$ such that each $B_n$ is $u_n$-bounded, which means that $B_n\subset\cup v$ for some finite
$v\subset u_n$. Substituting "$\gamma$" for "$\omega$" in the above sentence,
we obtain the definition of the property $\bigcup_{\mathrm{fin}}(\Gamma,\Omega)$,
which will be refered in this paper as the property of Scheepers. If, additionaly,
each $v_n$ in the definition of the Menger property contains only one element of $u_n$,
we obtain the definition of the Rothberger property. 

The following problem is open.
\begin{problem} (\cite[Issue 9, Problem 4.1]{SPM}, \cite[Problem 6.7]{Ts-new}). 
Is the property $\mathrm{Split}(\Lambda,\Lambda)$ preserved by unions of subsets of $\mathbb{R}$?
\end{problem}
We shall show that under additional strong set-theoretic assumption $\mathfrak u<\mathfrak g$  
every Lindel\"of paracompact topological space $X$ is Hurewicz provided it has the property
$\mathrm{Split}(\Lambda,\Lambda)$, which implies that positive answer onto the above problem is consistent.
In particular, this implies that under $\mathfrak u<\mathfrak g$ every Rothberger
space is Hurewicz. 
It is worth to mention here, that under CH there are  so called Luzin subsets of the
Baire space $\mathbb{N}^\omega$,
which have the Rothberger property but fail be Hurewicz, see \cite{JMSS} for details.
Therefore the statement "the family of Hurewicz subspaces  and
the family of subspaces with the property $\mathrm{Split}(\Lambda,\Lambda)$ of the real line coincide"
is independent of ZFC.

The reason why such the different results of these two parts are unified in one paper
is that both of them are proved with use of semifilters. 
\medskip

\large \centerline{\textbf{Semifilters}} \normalsize
\smallskip

 To begin with, let us recall from \cite{Va} the definition of the small cardinal
$\mathfrak g$.  Let $C$ be a countable set.
A family $\mathcal D\subset [C]^{\aleph_0}$ is said to be \emph{open}, if $X\in\mathcal D$ provided $X\subset^\ast Y$
for some $Y\in\mathcal D$ (here and subsequently $X\subset^\ast Y$ means that the complement 
$X\setminus Y$ is  finite, and $[A]^{\aleph_0}$ ($A^{<\aleph_0}$) denotes the set of all countable infinite 
(finite)
subsets of a set $A$). 
A family $\mathcal D$ is called \emph{groupwise dense}, if for every infinite 
collection $\Pi$ of finite parwise disjoint subsets of $C$
there exists an infinite  $\mathcal H\subset\Pi$ such that $\cup\mathcal H\in\mathcal D$.
By definition, $\mathfrak g$ equals to the smallest cardinality of a collection of
groupwise dense families with empty intersection.
Given arbitrary groupwise dense family $\mathcal D$, consider the family
$\mathcal F=\{C\setminus D:D\in\mathcal D\}\cup \mathfrak Fr(C)$, where $\mathfrak Fr(C)$ denotes the Fr\'echet filter
on $C$ consisting of cofinite subsets. From the above it follows that
$\mathcal F$ satisfies the following conditions:
\begin{itemize}
\item[$(1)$] $G\in\mathcal F$ provided $F\subset^\ast G$ for some $F\in\mathcal F$;
\item[$(2)$] every collection $\Pi$ of pairwise disjoint finite subsets 
of $C$ contains an infinite subset $\mathcal H$ such that $C\setminus \cup\mathcal H$
belongs to $\mathcal F$.
\end{itemize}
Following \cite{BZ}, we define a family $\mathcal F$ of infinite subsets of $C$
to be a \emph{semifilter}, if it satisfies the above mentioned condition $(1)$.

However, another approach to the definition of groupwise dense families 
is not the purpose of introduction of semifilters. Quite the contrary, semifilters
seem to constitute some rather interesting area of Set Theory, see \cite{BZ}.
In particular,  they  inherited many useful properties of filters, for example
the following classical theorem due to Talagrand holds, see \cite{Ta} or \cite{BZ}.
\begin{theorem} \label{Tal}
Let $\mathcal F$ be a semifilter on a countable set $C$. Then $\mathcal F$ fails to be  meager
if and only if it satisfies the above mentioned condition $(2)$.
\end{theorem} 
(Since every semifilter $\mathcal F$ on a countable set $C$ is a subset of the
powerset $\mathcal P(C)$, which can be identified with the product $\{0,1\}^C$,
we can speak about topological properties of semifilters. Since $C$ is countable,
$\mathcal P(C)$ and $[C]^{\aleph_0}$ are nothing else but homeomorphic copies
of the Cantor and Baire space respectively. For example, 
the base of the topology on $[C]^{\aleph_0}$ consists of subsets of the form $G(s,t)=\{A\in [C]^{\aleph_0}:A\cap s=t\}$,
where $s$ and $t$ are finite subset of $C$.)  

Theorem~\ref{Tal} implies the following characterization of groupwise dense families:
a family $\mathcal D\subset[C]^{\aleph_0}$ is groupwise dense if and only if 
the family $\{C\setminus D:D\in\mathcal D\}\cup\mathfrak Fr(C)$ is a nonmeager semifilter.
Therefore $\mathfrak g$ is equal to the smallest cardinality of a collection
$\mathsf F$ of nonmeager semifilters such that $\bigcap\mathsf F=\mathfrak Fr(C)$.
We shall prove a bit more.
\begin{observation}
The cardinal $\mathfrak g$ is equal to the smallest cardinality
of a family $\mathsf F$ of semifilters on a countable set $C$ such that $\bigcap \mathsf F$
is meager.
\end{observation}
\begin{proof}
Let $\mathsf F$ be a family of semifilters such that $\bigcap\mathsf F$ is meager.
The only thing to be proved is that $|\mathsf F|\geq\mathfrak g$.
For this aim let us fix a sequence $(I_n)_{n\in\omega}$ of pairwise disjoint
finite subsets of $C$ such that each member of $\bigcap\mathsf F$ meets
all but finitely many $I_n$. Without loss of generality, $\bigcup_{n\in\omega}I_n=C$.
For each $\mathcal F\in\mathsf F$ let us make the notation $\mathcal G_\mathcal{F}=
\{G\subset\omega:\bigcup_{n\in G}I_n\in\mathcal F\}$. Now, it sufficies to observe that
each $\mathcal G_\mathcal{F}$ is a nonmeager semifilter on $\omega$ and $\bigcap_{\mathcal F\in\mathsf F}\mathcal G_\mathcal{F}=\mathfrak Fr(\omega)$.
\end{proof}
The family of all semifilters on a set $C$ is evidently closed under taking unions and intersections
of arbitrary subfamilies. In addition to these operations there is another unary one.
Given any semifilter $\mathcal F$, let $\mathcal F^\perp=\{G\subset C:\mathcal \forall F\in\mathcal F(F\cap G\neq\emptyset)\}$.
(For a filter $\mathcal F$ the family $\mathcal F^\perp$ is nothing else but $\mathcal F^+$
in notations of C.Laflamme, see \cite{La}). It is clear that $\mathcal F^\perp$
is a semifilter too. In other words,
$\mathcal F^\perp=\mathcal{P}(C)\setminus\{C\setminus F:F\in\mathcal F\}$.
Consequently $(\mathcal F^\perp)^\perp=\mathcal F$ and $\mathcal F^\perp$ is comeager if and only if $\mathcal F$ is meager. Let us also observe
that $(\bigcap\mathsf F)^\perp=\bigcup_{\mathcal F\in\mathsf F}\mathcal F^\perp$ for arbitrary collection of
semifilters $\mathsf F$. Thus we obtain another characterization of the cardinal
$\mathfrak g$: it is the smallest size of a family $\mathsf F$ of non comeager semifilters on a countable set $C$
such that $\bigcup\mathsf F$ is comeager. In what follows we shall simply write
$\mathfrak Fr$ in place of $\mathfrak Fr(\omega)$.

Next, similarly to \cite{BM} for every semifilter $\mathcal F$ on $\omega  $ we shall define 
a cardinal characteristic $\mathfrak b(\mathcal F)$. Its definition involves some special relation
$\leq_\mathcal F$ on $\mathbb N^\omega$:
\[  (x_n)_{n\in\omega}\leq_\mathcal F(y_n)_{n\in\omega} \mbox{ iff } \{n\in\omega:x_n\leq y_n\}\in\mathcal F.           \]  
Now, $\mathfrak b(\mathcal F) $ stands for the smallest size of unbounded subset of $\mathbb N^\omega$ with respect
to $\leq_\mathcal F$. When $\mathcal F=\mathfrak Fr$, then $\leq_{\mathcal F}$ is nothing
else but the well-known eventual dominance preorder $\leq^\ast$.
For example, $\mathfrak b(\mathfrak Fr^\perp)=\mathfrak d$ and $\mathfrak b(\mathfrak Fr)=\mathfrak b$. 
Almost literal repetition of the proof of Theorem~16 from \cite{BM} gives us the following
\begin{proposition} \label{fromBM}
$\mathfrak b(\mathcal F)\geq\mathfrak g$ for each nonmeager semifilter $\mathcal F$ on $\omega$.
\end{proposition}
Next, in what follows we shall intensively use set-valued maps. 
By a \emph{set-valued map} $\Phi$ from a set $X$
into a set $Y$ we understand a map from $X$ into $\mathcal{P}(Y)$ and write
$\Phi:X\Rightarrow Y$ (here $\mathcal P(Y)$ denotes the set of all subsets of $Y$). For a subset $A$ of $X$ we put $\Phi(A)=\bigcup_{x\in A}\Phi(x)\subset Y.$
 When the sets $X$ and $Y$ are endowed with some topologies,
it is interesting to consider set-valued maps with certain topological properties.
The set-valued map $\Phi$ between topologival spaces $X$ and $Y$ is said to be
\begin{itemize}
\item \emph{compact-valued}, if $\Phi(x)$ is compact for every $x\in X$;
\item\emph{upper semicontinuous}, if for every open subset $V$ of $Y$ the set
$\Phi_\subset^{-1}(V)=\{x\in X:\Phi(x)\subset V\}$ is open in $X$.
\end{itemize}
\begin{lemma} \label{mf}
Let $\Phi:X\Rightarrow Y$ be a compact-valued upper semicontinuous map
between topological spaces $X$ and $Y$ such that $\Phi(X)=Y$. Then $Y$ is Menger (Hurewicz) provided so is $X$.
\end{lemma}
\begin{proof}
Let  us fix arbitrary sequence $(w_n)_{n\in \omega}$ of open covers of $Y$.
For every $n\in\omega$ consider the family $u_n=\{\Phi_\subset^{-1}(\cup v):v\in [w_n]^{<\aleph_0}\}$.
Since $\Phi$ is upper semicontinuous and compact-valued, each $u_n$ is an open cover of $X$.
The Menger property of $X$ implies the existence of a sequence $(c_n)_{n\in\omega}$,
where each $c_n$ is a finite subset of $u_n$, such that $\{\cup c_n:n\in\omega\}$
is a ($\gamma$-) cover of $X$. From the above it follows that for every $n\in\omega$ we can find
a finite subset $v_n$ of $w_n$ with $\Phi(\cup c_n)\subset \cup v_n$. Therefore
$\{\cup v_n:n\in\omega\}$ is  a ($\gamma$-) cover of $Y$,
consequently $Y$ is Menger (Hurewicz).
\end{proof}

The main idea of this paper is to assign to a topological space $X$ the collection $\mathsf U(X)=\{\mathcal U(u,X):u\in\Lambda_\omega(X)\}$
of semifilters on countable sets, where $\Lambda_\omega(X)$  denotes the family
of all countable large open covers of $X$ 
 and $\mathcal U(u,X)$ is the smallest semifilter
on $u$ containing the family $\{I(x,u,X)=\{U\in u:x\in U\}:x\in X\}$.
It is clear that $\mathcal U(u,X)$ can be represented in the form
$\bigcup_{v\in [u]^{<\aleph_0}}\bigcup_{x\in X}\uparrow_v I(x,u,X)$, where
for a subsets $A$ and $B$ of a set $Z$ we denote by $\uparrow_A B$ the family
$\{C\subset Z:C\supset B\setminus A\}$. When $A=\emptyset$, we shall simply write
$\uparrow$ in place of $\uparrow_A$. When $X$ (and $u$) are clear from the context,
we shall write $\mathcal U(u)$ and $I(x,u)$ ($I(x)$) instead of $\mathcal U(u,X)$ and $I(x,u,X)$.

We are in a position now to present a characterization of the properties of
Menger and Hurewicz in terms of topological properties of semifilters,
which implies the mentioned in Introduction results.

\begin{theorem} \label{ch0}
Let $X$ be a Lindel\"of topological space. Then $X$ is Menger (Hurewicz)
if and only if so is each $\mathcal U(u)\in\mathsf U(X)$.
Moreover, if $X$ is paracompact, then it is Hurewicz provided each semifilter
$\mathcal U(u)\in \mathsf U(X)$ is meager.
\end{theorem}
\noindent \textbf{Remark}. 1. Every Hurewicz semifilter on a countable set $C$ is meager.
Indeed, \cite[Theorem 5.7]{JMSS} implies that each Hurewicz semifilter $\mathcal F$
on $C$ is contained in a $\sigma$-compact subset of $[C]^{\aleph_0}$,
and each $\sigma$-compact subset of the Baire space is meager.

2. The ``meager'' part of the above characterization of the Hurewicz property was 
recently obtained by B.~Tsaban \cite{sl} in case of a zero-dimensional metrizable
space $X$.
\hfill $\Box$
\medskip

We shall divide the proof of Theorem~\ref{ch0} into a sequence of lemmas.
\begin{lemma}\label{crucial}
Let $X$ be a topological space and $u\in\Lambda_\omega(X)$. 
Then the  set-valued map $\Phi:X\Rightarrow\mathcal P(\omega)$,
$\Phi: x\mapsto \uparrow I(x)$,
is compact-valued and upper semicontinuous.
\end{lemma}
\begin{proof}
It is clear that $\Phi$ 
is compact-valued, because $\Phi(x)=\uparrow I(x)$ is closed and precompact subspace of
$\mathcal P(u)$. 
Let us show that $\Phi$ is upper semicontinuous. For this aim let us consider arbitrary $x\in X$
and an open subset $G$ of $\mathcal P(u)$ containing $\Phi(x)$.
 For every $v\in\Phi(x)$ we can find   
$s_v\in [u]^{<\aleph_0}$ such that $G(s_v,s_v\cap v)\subset G$. Since $\Phi(x)$ is compact,
we can find a finite family $\mathbf v\subset\Phi(x)$ such that $\Phi(x)\subset\bigcup\{G(s_v,s_v\cap  v):v\in\mathbf v\}$.
Put $s=\cup\{s_v:v\in\mathbf v\}$,  $c=s\cap I(x)$ and $U=\cap c$.  
It is clear that $x\in U$ and $U$ is  open. Thus the upper semicontinuity of $\Phi$ 
will be proven as soon as we show that $\Phi(U)\subset G$. For this purpose let us fix 
arbitrary $x_1\in U$ and observe that $I(x_1)\cap s\supset c=I(x)\cap s$,
consequently $\Phi(x_1)\subset\bigcup\{G(s,v\cap s):v\in\Phi(x)\}=G$,
and finally $\Phi(U)\subset G$, which implies upper semicontinuity of $\Phi$.
\end{proof}

\begin{corollary}\label{ccrucial}
Let $X$ be a Menger (Hurewicz) topological space and $u\in\Lambda_\omega(X)$. Then the semifilter
$\mathcal U(u,X)$ is Menger (Hurewicz). 
\end{corollary}
\begin{proof}
Given any $v\in [u]^{<\aleph_0}$, consider the set-valued map $\Phi_v:X\Rightarrow\mathcal P(u)$,
$\Phi_v:x\mapsto\uparrow_v I(x,u)$. Let us observe, that 
$\Phi_v$ is a composition $\Psi_2\circ\Psi_1$, where $\Psi_1:X\Rightarrow\mathcal P(u\setminus v)$,
$\Psi_1:x\mapsto\uparrow I(x,u\setminus v)$, and $\Psi_2:\mathcal P(u\setminus v)\Rightarrow\mathcal P(u)$,
$\Psi_2:w\mapsto \uparrow w$. It is clear that $\Psi_2$ is compact-valued upper
semicontinuous, while $\Psi_1$ is so by Lemma~\ref{crucial}. Now,
Lemma~\ref{mf} implies that $\bigcup_{x\in X}\uparrow_v I(x,u)=\Phi_v(X)$
is Menger (Hurewicz). Since the property of Menger (Hurewicz) is preserved by
countable unions, the semifilter $\mathcal U(u)=\bigcup_{v\in [u]^{<\aleph_0}}\Phi_v(X)$
is Menger (Hurewicz).
\end{proof}

\begin{lemma} \label{crucial1}
Let $X$ be a Lindel\"of topological space which fails to be Menger \\ 
(Hurewicz).
 Then there exists $u\in\Lambda_\omega(X)$ such that
the semifilter $\mathcal U(u)$ is not
Menger (Hurewicz).
\end{lemma}
\begin{proof}
Assuming that $X$ is not Menger (Hurewicz), we can find a sequence $(u_n)_{n\in\omega}$
of countable open large covers of $X$ such that there is no sequence
$(v_n)_{n\in\omega} $ such that each $v_n$ is a finite subset of $u_n$
and the family $\{\cup v_n:n\in\omega\}$ is a ($\gamma$-)cover of $X$.
Let us  denote by $u$ the union $\cup\{u_n:n\in\omega\}$.

We claim that the semifilter $\mathcal U(u)$ 
is not Menger (Hurewicz). Indeed, consider the 
sequence $(o_n)_{n\in\omega}$ of countable families of open subsets of $\mathcal P(u)$,
 where $o_n=\{\{w\in \mathcal P(u):U\in w\}:U\in u_n\}$.
Since each $u_n$ is a large cover of $X$, every $o_n$ covers $\mathcal U(u)$.
It sufficies to show that there is no sequence $(c_n)_{n\in\omega}$ such that
every $c_n$ is a finite subset of $o_n$ and $\{\cup c_n:n\in\omega\}$ is a large ($\gamma$-) cover of $\mathcal U(u)$, see \cite[Corollary 5]{Sch}.
Assume, to the contrary, that such the sequence $(c_n)_{n\in\omega}$ exists. Then for every
$n\in\omega$ we can find a finite subset $v_n$ of $u_n$ such that $c_n=\{\{w\in\mathcal P(u):w\ni U\}:U\in v_n\}$.
For every $w\in\mathcal U(u)$ set $J_w=\{n\in\omega:w\in\cup c_n\}=\{n\in\omega:w\cap v_n\neq\emptyset\}$.
From the above it follows that $\mathcal J=\{J_w:w\in\mathcal U(u)\}$ consists of infinite (cofinite) subsets
of $\omega$.
From the above it follows that the family $\{J_{I(x,u)}:x\in X\}$
consists of infinite (cofinite) subsets of $\omega$ too. But
\[ 
J_{I(x,u)}=\{n\in\omega:I(x,u)\cap v_n\neq\emptyset\}=\{n\in\omega:x\in\cup v_n\},\] 
consequently  $\{\cup v_n:n\in\omega\}$ is a ($\gamma$-) cover of $X$,
which contradicts our choice of the sequence $(u_n)_{n\in\omega}$.
\end{proof}

Let $u$ be a family of a set $X$ and $B\subset X$. From now on
$\mathcal St(B,u)$ denotes the set $\cup\{U\in u:U\cap B\neq\emptyset\}$.

 \begin{lemma}
\label{crucial2}
Let  $X$ be 
 a paracompact Lindel\"of topological space. 
 Then $X$ is Hurewicz provided
 each semifilter $\mathcal U(u)\in\mathsf U(X)$
is meager.
\end{lemma}
\begin{proof}
Assuming that $X$ is not Hurewicz, we shall show that $X$ possess
a countable large open cover $u$ such that the semifilter $\mathcal U(u)$ is not meager. 
Let $(u_n)_{n\in\omega}$ be a sequence
of open covers of $X$ such that $\{\cup v_n:n\in\omega\}$ is a $\gamma$-cover of $X$ for no sequence
$(v_n)_{n\in\omega} $ such that each $v_n$ is a finite subcolection of $u_n$.
Now, it is a simple exercise to construct a sequence $(w_n)_{n\in\omega}$ of open covers of $X$,
where $w_n=\{U_{n,k}:k\in\mathbb{N}\}$ is a refinement of $u_n$, such that  
$U_{n_2,k}\subset\bigcup\{U_{n_1,l}:l\leq k\}$ for all $n_2\geq n_1$
 and $\mathcal St(B,w_n)$ is $w_n$-bounded for every $w_n$-bounded subset $B$ of $X$.
From the above it follows that for every $n\in\omega$ there exists a sequence
 $(p(n,k))_{k\in\omega} $ of natural numbers such that $\cup\{U_{m,i}:i\leq k, n\leq m\}\bigcap\cup\{U_{n,i}:i\geq p(n,k)\}=\emptyset$.
Without loss of generality, $u_n=w_n$. 

 Denote by $u$ the union $\cup\{u_n:n\in\omega\}$. 
Let $\nu:u\rightarrow\omega$, $\nu:U_{n,k}\mapsto m_{n,k}$ be a bijective enumeration of $u$.
Let us write $\omega$ in the form $\omega=\sqcup_{n\in\omega}\Omega_n$, where $\Omega_n=\{m_{n,k}:k\in\mathbb{N}\}$.
We claim that the semifilter $\mathcal U(u)$ fails to be meager. 
For this aim we shall show that the image $\mathcal F=\nu(\mathcal U(u))$ 
of $\mathcal U(u)$ is not meager in $\mathcal P(\omega)$.
Otherwise, by Theorem~\ref{Tal} there exists a sequence $(m_l)_{l\in\omega}$ of natural numbers 
such that every $F\in\mathcal F$ (and, in particular,  $\nu(I(x,u))$ for every $x\in X$) meets all but finitely many half-intervals
 $[m_l,m_{l+1})$.
  Passing to a subsequence, if necessary, we may assume that
$m_{l+1}>\max\{m_{n_1,p(n_2,k)}:m_{n_2,k}\leq m_l, [0,m_l]\cap\Omega_{n_1}\neq\emptyset\}$.
Consider a  function $\varphi:\omega\rightarrow\omega$
such that $\varphi^{-1}(l)=[m_l,m_{l+1})$ for all $l\in\omega$ and denote by
 $B_l$ the union $\cup\{\nu^{-1}(m):m\in\varphi^{-1}(l)\}$.  
Then for every $x\in X$ there exists $n\in\omega$ such that $(\varphi\circ\nu)(I(x,u))\supset [n,+\infty)$.
From the above it follows that 
$$X=\bigcup_{x\in X}\bigcap_{m\in \nu(I(x,u))}\nu^{-1}(m)\subset\bigcup_{x\in X}\bigcap_{l\in(\varphi\circ\nu)(I(x,u))}\bigcup_{m\in\varphi^{-1}(l)}\nu^{-1}(m)\subset\bigcup_{n\in\omega}\cap_{i\geq n}B_i.$$
A crucial observation here is that the intersection $\bigcap_{l\in A}B_l$ is $u_n$-bounded
for every $n\in\omega$ and all infinite subsets $A$ of $\omega$. 
Before proving this observation,
let us note, that we can reduce ourselves to subsets $A$ such that $|l_1-l_2|>1$
for all $l_1,l_2\in A$. 
Given arbitrary $l\in\omega$,
 denote by 
$K_{n,l}$ the intersection $\Omega_n\cap[m_l,m_{l+1})$.
 Equipped with these notations, we can write
$$\bigcap_{l\in A}B_l=\bigcap_{l\in A}\bigcup_{n\in\omega}\bigcup_{m\in[m_l,m_{l+1})\cap\Omega_n}\nu^{-1}(m)=\bigcap_{l\in A}\bigcup_{n\in\omega}\bigcup_{m\in K_{n,l}}\nu^{-1}(m).$$
Set $B_{n,l}=\cup\{\nu^{-1}(m):m\in K_{n,l}\}$. Then
$\bigcap_{l\in A}B_l=\bigcap_{l\in A}\cup_{n\in\omega}B_{n,l}=\\
=\bigcup_{z\in\mathbb{N}^A}\cap_{l\in A}B_{z(l),l}.$
By our choice of a sequence $(m_l)_{l\in\omega}$ and the subset $A$ of $\omega$ 
 we have  
$B_{z(l_1),l_1}\cap B_{z(l_2),l_2}=\emptyset$ provided $z(l_2)\leq z(l_1)$ for some $l_1,l_2\in A$
with $l_2>l_1$. Consequently
$$\bigcap_{l\in A}B_l=\bigcup_{z\in\mathbb{N}^{\uparrow A}}\cap_{l\in A}B_{z(l),l}=\bigcap_{l\in A}\bigcup_{n\geq |A\cap [0,l)|}B_{n,l}, $$
where $\mathbb{N}^{\uparrow A}=\{z\in\mathbb{N}^A:z(l_2)>z(l_1)\mbox{ for all }l_2>l_1\}$. 
Since the union \\ 
$\bigcup_{n\geq |A\cap [0,l)|}B_{n,l}$
is $u_{|A\cap [0,l)|}$-bounded and $|A\cap [0,l)|\rightarrow +\infty$, $l\rightarrow +\infty$,
the above intersection is $u_n$-bounded for all $n\in\omega$ (recall, that
 each $u_{n+1}$-bounded subset of $X$ is $u_n$-bounded). Therefore there exists a sequence
$(v_n)_{n\in\omega}$, where $v_n$ is a finite subset of $u_n$, such that
 $\cup v_n\supset\bigcup_{k\leq n}\cap_{i\geq k}B_i$,
 consequently the family $\{\cup v_n:n\in\omega\}$ is a $\gamma$-cover of $X$,
 which contradicts our choice of the sequence
$(u_n)_{n\in\omega}$.
\end{proof}
\medskip

\textbf{Proof of Theorem~\ref{ch0}}.
Follows from Lemmas~\ref{mf}, \ref{crucial1}, \ref{crucial2}, Corollary~\ref{ccrucial},
and the remark after the formulation of Theorem~\ref{ch0}.
\hfill $\Box$
\medskip

The following statement is of great importance in evaluation 
of additivity of the family of subspaces with the Menger
property of a topological space $X$.
\begin{proposition} \label{sm:sem}
Every comeager semifilter $\mathcal F$ on $\omega$ is not Menger.
\end{proposition} 
\begin{proof}
Since $\mathcal F$ is comeager in the space $[\omega]^{\aleph_0}$, which is homeomorhic to the
Baire space $\mathbb N^\omega$, there exists a dense $G_\delta $ subset $G$ of $[\omega]^{\aleph_0}$ 
such that $G\subset\mathcal F$. Thus $G$ is an analytic and not $\sigma$-compact subset of
$[\omega]^{\aleph_0}$, consequently it contains a closed in $[\omega]^{\aleph_0}$ subset $D$
homeomorphic to $\mathbb N^\omega$, see \cite[Theorem 29.3]{Ke}. But $\mathbb N^\omega$ simply fails to be Menger,
consequently so is $\mathcal F$, a contradiction. 
\end{proof}

Theorem~\ref{ch0} and Proposition~\ref{sm:sem} enable us to introduce a new class of topological spaces.
A topological space $X$ is defined to be \emph{almost Menger}, if semifilter $\mathcal U(u)$
is not comeager for every $u\in\Lambda_\omega(X)$. Theorem~\ref{ch0}
implies that every Lindel\"of Menger space is almost Menger.
\begin{problem} \label{am=m}
Is every (metrizable separable) Lindel\"of almost Menger space\\
 Menger?
\end{problem} 
\medskip

Sometimes  it is more convenient to use some modification of Theorem~\ref{ch0}.
Let $X\subset Y$ and $u=(U_n)_{n\in\omega}$ be a sequence of subsets of $Y$.
For every $x\in X$ let  $I_s(x,u,X)=\{n\in\omega:x\in U_n\}$. If every
$I_s(x,u,X)$ is infinite, then we shall denote by $\mathcal U_s(u,X)$
the smallest semifilter on $\omega$ containing all $I_s(x,u,X)$ (the letter $s$
comes from "sequence"). In what follows we shall denote by $\Lambda_s(X)$
the set of all sequences $u=(U_n)_{n\in\omega}$ of open nonempty subsets of a topological space $X$
such that all $I_s(x,u X)$ are infinite. 
 Again, we shall often simplify these notations by
writing $\mathcal U_s(u)$ and $I_s(x,u)$ or $I_s(x)$ in place of $\mathcal U_s(u,X)$
and $I_s(x,u,X)$. 

\begin{theorem} \label{ch1}
Let $X$ be a Lindel\"of topological space. Then $X$ is Menger (Hurewicz) if and only if
for every sequence $u=(U_n)_{n\in\omega}\in\Lambda_s(X)$  the semifilter $\mathcal U_s(x,u)$ is Menger (Hurewicz).
In addition, if $X$ is paracompact, then it is Hurewicz provided $\mathcal U_s(u)$ is meager
for every $u\in\Lambda_s(X)$.
\end{theorem}
\begin{proof}
Assuming that $X$ is Menger (Hurewicz), let us fix any sequence $u=(U_n)_{n\in\omega}\in\Lambda_s(X)$.
and denote by $Y$ the product $X\times\mathbb{N}$. The space $Y$ is Menger (Hurewicz)
being a countable union of its Menger (Hurewicz) subspaces. Consider the cover
$w=\{W_n:n\in\omega\}$ of $Y$, where $W_n=U_n\times\{1,\ldots,n\}$, and observe that
$w\in\Lambda_\omega(Y)$. Applying Theorem~\ref{ch0}, we conclude that $\mathcal U(w)$ is Menger
(Hurewicz) subspace of $\mathcal P(w)$. Now, it sufficies to observe that
$\mathcal U_s(u)$ is a continuous image of $\mathcal U(w)$ under the map $f:w\rightarrow\omega,$
$f:W_n\mapsto n$.

Next, assume that the semifilter $\mathcal U_s(u)$ is Menger (resp. Hurewicz, meager)
for all $u\in\Lambda_s(X)$ and fix any $w\in\Lambda_w(X)$.
Let $w=\{W_n:n\in\omega\}$ be a bijective enumeration of $w$ and $u$ be a sequence
$(W_n)_{n\in\omega}$. Then $f(\mathcal U_s(u))=\mathcal U(w)$, where $f:n\mapsto W_n$ is a bijection.
Therefore the semifilter $\mathcal U(w)$ is Menger (resp. Hurewicz, meager). Now, it sufficies to apply
Theorem~\ref{ch0}.
\end{proof}
 
\large \centerline{\textbf{Additivity of the Menger property}} \normalsize
\smallskip

As we have already said in Introduction, one of the main result of this paper is the following
\begin{theorem} \label{main}
Let $X$  a hereditarily Lindel\"of
space. Then \\
$\mathrm{add}(M(X))\geq\mathrm{add}(M(\mathbb{N}^\omega))\geq\mathfrak g$.
\end{theorem}
\begin{proof}
Let $\mathcal Y$ be a subfamily of $M(X)$ of size $|\mathcal Y|<\mathrm{add}(M(\mathbb{N}^\omega))$
and $u=(U_n)_{n\in\omega}\in\Lambda_s(\cup\mathcal Y)$.  Then
the semifilter $\mathcal U_s(u,\cup\mathcal Y)$ is equal to the union \\
$\bigcup_{Y\in\mathcal Y}\mathcal U_s(u,Y)$.
Theorem~\ref{ch1} implies that every $\mathcal U_s(u,Y)$ is Menger, consequently
so is their union $\mathcal U_s(u,\cup\mathcal Y)$ by our choice of $\mathcal Y$.
Applying Theorem~\ref{ch1} once again, we conclude that $\cup\mathcal Y$
is Menger, 
 which implies the inequality
$\mathrm{add}(M(X))\geq\mathrm{add}(M(\mathbb{N}^\omega))$.

Next, we shall show that $\mathrm{add}(M(\mathbb{N}^\omega))\geq\mathfrak g$.
Let $(w_n)_{n\in\omega}$ be a sequence of open covers of $\bigcup\mathcal Y$.
Since $X$ is hereditarily Lindel\"of, we can assume that every $w_n$ is a countable cover of $\bigcup \mathcal Y$
of the form 
$w_n=\{W_{n,k}:k\in\omega\}$. For every $Y\in\mathcal Y$ let us find a sequence
$(k_n(Y))_{n\in\omega}$ of natural numbers such that the sequence $u_Y=(B_n(Y))_{n\in\omega}$, where $B_n(Y)=\bigcup_{k\leq k_n(Y)}W_{n,k}$,
belongs to $\Lambda_s(Y)$. Now, Theorem~\ref{ch1} and Proposition~\ref{sm:sem}  imply that the  semifilter
$\mathcal U_s(u_Y,Y)$ fails to be comeager. Since $|\mathcal Y|<\mathfrak g$,  the semifilter $\mathcal F=\bigcup_{Y\in\mathcal Y}\mathcal U_s(u_Y,Y)$
is not comeager too, consequently
 $\mathcal F^\perp$ fails to be meager. Using $|\mathcal Y|<\mathfrak g\leq\mathfrak b(\mathcal F^\perp)$, we can find
a sequence $(k_n)_{n\in\omega}$ such that $(k_n(Y))_{n\in\omega}\leq_{\mathcal F^\perp}(k_n)_{n\in\omega}$
for every $Y\in\mathcal Y$. Let us make the following notation: $B_n=\bigcup\{U_{n,k}:k\leq k_n\}$.
We claim that $\{B_n:n\in\omega\}$ is a cover of $\cup\mathcal Y$. Indeed, let us fix arbitrary $Y\in\mathcal Y$
and $y\in Y$. Since
$(k_n(Y))_{n\in\omega}\leq_{\mathcal F^\perp}(n_k)_{k\in\omega}$, the set $A=\{n\in\omega:k_n\geq k_n(Y)\}$
belongs to $\mathcal F^\perp$ and thus there exists $m\in A\cap I_s(y,u_Y,Y) $. It sufficies to observe that $y\in B_m$.
\end{proof}

\begin{problem}
Is the equation $\mathrm{add}(M(X))=\mathrm{add}(M(\mathbb{N}^\omega))$ true for every hereditarily Lindel\"of topological space? 
\end{problem}
\medskip

\large \centerline{\textbf{Additivity of the property $\mathrm{Split}(\Lambda,\Lambda)$}} \normalsize
\smallskip

Throughout this paragraph, which is devoted to the property $\mathrm{Split}(\Lambda,\Lambda)$,
 every topological space is hereditarily Lindel\"of.
Since every large open cover of a  space $X$ contains a countable large subcover (see 
\cite[Proposition 1.1]{Ts3}),
we can restrict ourselves to countable ones. 
\begin{theorem}\label{split}
Under $\mathfrak u<\mathfrak g$ every paracompact space $X$ with the property
$\mathrm{Split}(\Lambda,\Lambda)$ is Hurewicz.
\end{theorem}

In our proof of Theorem~\ref{split}  we shall use the following straightforward
consequence  of a fundamental result
of C.Laflamme. A semifilter $\mathcal F$ on a countable set $C$ is said to be
\emph{bi-Baire}, if it is neither meager nor comeager.
\begin{theorem}(\cite[Theorem 9.22]{Bl}, \cite{Laf1}) \label{lll}
 Let $C$ be a countable set and $\mathcal F$ be a semifilter on $C$.
If $\mathcal F$ is comeager ($\mathfrak u<\mathfrak g$ and $\mathcal F$ is bi-Baire), then there exists a sequence $(K_n)_{n\in\omega}$
of pairwise disjoint finite subsets of $C$ such that 
the set $\mathcal U=\{\{n\in\omega:F\cap K_n\neq\emptyset\}:F\in\mathcal F\}$ 
equals to $\mathfrak Fr^\perp$ (is an ultrafilter on
$\omega$).
\end{theorem}
\textbf{Remark.} Let us observe, that if a sequence $(K_n)_{n\in\omega}$
is such as in  Theorem~\ref{lll}, then for every
increasing sequence $(m_n)_{n\in\omega}$ of natural numbers the sequence $(K_n'=\bigcup_{m\in[m_n,m_{n+1})}K_m)_{n\in\omega}$
satisfies the condition of this Theorem too.
\hfill $\Box$.
\medskip

\textbf{Proof of Theorem~\ref{split}.}
In light of Corollary 29 from \cite{Sch} asserting that each Hurewicz space
has the property $\mathrm{Split}(\Lambda,\Lambda)$, the only step to be proven is the inverse implication under
$\mathfrak u<\mathfrak g$.  Supppose that the paracompact space $X$ 
is not Hurewicz. Then Theorem~\ref{ch0} supplies us with a cover $u\in\Lambda_\omega(X)$
such that the semifilter $\mathcal U(u)$ is not meager. Therefore there exists a finite
subset $v$ of $u$ such that no finite subset $v_1$ of $u\setminus
v$ is a cover of $X$, because otherwise we can simply construct by induction
a sequence $(v_n)_{n\in\omega}$ of pairwise disjoint finite subsets of $u$
such that each $w\in\mathcal U(u)$ meets all but finitely many $v_n$,
and thus $\mathcal U(u)$ is meager by Theorem~\ref{Tal}. 
Two cases are possible.

1.  $\mathcal U(u)$ is comeager. Then we can find a sequence $(v_n)_{n\in\omega}$ of 
finite subsets of $u$ such as in Theorem~\ref{lll}. 
Let $n_0\in\omega$ be such that no finite subset of $\cup_{n\geq n_0}v_n$
covers $X$. Since every $w\in\mathcal U(u)$ meets $u_0=\cup_{n\geq n_0}v_n$,
we conclude that $u_0\in\Lambda_\omega(X)$. From the above it follows that
there exists an increasing sequence $(m_n)_{n\in\omega}$ of natural numbers
such that $m_0\geq n_0$ and $U_{n_1}\neq U_{n_2}$ for all $n_1\neq n_2$, where   $U_n=\bigcup_{m\in [m_n,m_{n+1})}\cup v_m,$ $n\in\omega$.
Let us denote by $v_n'$ the union $\cup_{m\in[m_n,m_{n+1})}v_m$ and 
observe that the sequence $(v'_n)_{n\in\omega}$ satisfies the condition of Theorem~\ref{lll}.

It is clear that
the family $u'=\{U_n=\cup v'_n:n\in\omega\}$ is a large cover of $X$. We claim
that $u'$  is not
splittable. Assuming the converse, we would find two disjoint infinite subsets
$A$ and $B$ of $\omega$  such that both of $u_A=\{U_n:n\in A\}$ and $u_B=\{U_n:n\in B\}$
are large covers of $X$. Since the sequence $(v'_n)_{n\in\omega}$ satisfies
the condition of Theorem~\ref{lll}, there exists $w\in \mathcal U(u)$
such that $\{n\in\omega:w\cap v_n'\}=A$. By the definition of the semifilter $\mathcal U(u)$,
$I(x,u)\subset^\ast w$ for some $x\in X$. Therefore 
the set $\{n\in B:x\in U_n\}$ is a subset of a finite set
$\{n\in B:(I(x,u)\setminus w)\cap v_n'\neq\emptyset\}$, and thus $u_B$
is not a large cover of $X$, a contradiction.

2. $\mathcal U(u)$ is not comeager. Then  the same argument     as in the first
case implies that there exists a sequence $(v_n)_{n\in\omega}$ of finite subsets
of $u$ with the following properties:
\begin{itemize}
\item[$1)$] it saatisfies the conditions of Theorem~\ref{lll};
\item[$2)$] no finite subset of $\cup_{n\in\omega}v_n$ covers $X$;
\item[$3)$] $\cup v_n\neq\cup v_m$, $n\neq m$, and the family $u'=\{\cup v_n:n\in\omega\}$
is a large cover of $X$.
\end{itemize}
We claim that $u'$ is not splittable. Assume, to the contrary, that
there are infinite disjoint subsets $A$ and $B$ of $\omega$
such that  $u_A,u_B\in\Lambda_\omega(X)$, where $u_A$ and $u_B$ are defined as above.
Enlarging $A$, if necessary, we can additionaly assume that $A\cup B=\omega$.
Since the family $\mathcal F=\{\{n\in\omega:w\cap v_n\neq\emptyset\}:w\in\mathcal U(u)\}$ is an ultrafilter,
either $A\in\mathcal F$ or $B\in\mathcal F$. Without loss of generality, $A\in\mathcal F$,
which means that there exists $x\in X$ such that $\{n\in\omega:I(x,u)\cap v_n\neq\emptyset\}\subset^\ast A$,
and thus $\{n\in B:I(x,u)\cap v_n\neq\emptyset\}=\{n\in B:x\in \cup v_n\}$ is finite,
a contradiction. 
\hfill $\Box$
\medskip

\large \centerline{\textbf{Another applications of Theorem~\ref{ch0}.}} \normalsize
\smallskip

Here we shall show that we can not  restrict ourselves to $\omega$-covers  in Theorem
\ref{ch0}, and thus the family $\mathsf U(X)$ of semifilters can not be
reduced to the family $\{\mathcal U(u):u\in\Omega(X)\}$ of filters corresponding to $\omega$-covers of a space
$X$. Thus Theorems~\ref{ch0} and \ref{ch1} are "purely semifilter statements".

\begin{proposition} \label{sa:re}
Let $u$ be an $\omega$-cover of $\mathbb{N}^\omega$. Then $\mathcal U(u)$ is meager.
Moreover, the smallest filter $\mathcal V$ containing $\mathcal U(u)$ is meager.
\end{proposition}
\begin{proof}
Since $u$ is an $\omega$-cover, the semifilter $\mathcal U(u)$ is centered and the
filter $\mathcal V$ is free. As it was shown in the proof of Corollary~\ref{ccrucial},
there exists a sequence $(\Phi_n)_{n\in\omega}$ of compact-valued upper semicontinuos
multifunctions from $\mathbb{N}^\omega$ into $\mathcal P(u)$ such that $\mathcal U(u)=\bigcup_{n\in\omega}\Phi_n(\mathbb{N}^\omega)$.
Each $\Phi_n(\mathbb{N}^\omega)$ is an image of $\mathbb{N}^\omega$ under a compact-valued upper
semicontinuous set-valued map ($\equiv$ $\Phi_n(\mathbb{N}^\omega)$ is $K$-analytic).
Since every $K$-analytic metrizable space $X$ is analytic (see \cite{JR}),
so is the semifilter $\mathcal U(u)$ being a countable union of analytic spaces,
see \cite[25.A]{Ke}. Let us note, that $\mathcal V=\bigcup_{n\in\omega}\mathcal U_n$, where
$\mathcal U_0=\mathcal U(u)$ and $\mathcal U_{n+1}$ is a continuous image of $\mathcal U_n^2$
under the map $(U_1,U_2)\mapsto U_1\cap U_2$.
From the above it follows that $\mathcal V$ is analytic too, consequently by
\cite[Theorem 1, p. 30]{To} it has the Baire property in $\mathcal P(u)$,
and thus is meager by \cite[Theorem 1, p. 32]{To}.
\end{proof}
Reformulating the above proposition in other terms, we obtain the subsequent result proved in
\cite{Sa2}.
\begin{theorem}\label{sasa} 
Every  $\omega$-cover of $\mathbb{N}^\omega$ is $\omega$-groupable.
\end{theorem}
Proposition~\ref{sm:sem} and Theorem~\ref{lll} imply the subsequent
\begin{corollary}
Under $\mathfrak u<\mathfrak g$ every Menger space is Scheepers ($\equiv$ has the property
$\bigcup_{\mathrm{fin}}(\Gamma,\Omega)$).
\end{corollary}
\begin{proof}
Let $(u_n)_{n\in\omega}$ be a sequence of open covers of $X$ such that $u_{n+1}$
is a refinement of $u_n$ for all $n\in\omega$. Since $X$ is Menger,
there exists a sequence $(v_n)_{n\in\omega}$ such that each $v_n$ is a finite subset of $u_n$
and $w=(\cup v_n)_{n\in\omega}$ belongs to $\Lambda_s(X)$. Applying Theorem~\ref{ch1} and Proposition~\ref{sm:sem}, we conclude that
the semifilter $\mathcal U_s(w)$ is not comeager. 

If $\mathcal U_s(w)$ is meager, then Theorem~\ref{Tal} gives us an increasing sequence 
$(m_n)_{n\in\omega}$ of natural numbers such that each $A\in\mathcal U_s(w)$ meets all but finitely many
half-intervals $[m_n,m_{n+1})$. Let $B_n=\bigcup_{m\in[m_n,m_{n+1})}\cup v_m$, $n\in\omega$.
From the above it follows that each $B_n$ is $u_n$-bounded and the family
$\{B_n:n\in\omega\}$ is a $\gamma$-cover of $X$.

If $\mathcal U_s(w)$ is bi-Baire,  then by Theorem~\ref{lll} there exists a sequence
$(K_n)_{n\in\omega}$ of finite subsets of $\omega$ such that the family $\mathcal F=\{\{n\in\omega:A\cap K_n\neq\emptyset\}:A\in\mathcal U_s(w)\}$
is an ultrafilter on $\omega$. Let $(m_n)_{n\in\omega}$  be an increasing sequence of natural numbers
with the property $\min\cup\{K_m:m\in[m_{n+2},m_{n+3})\}>\max\cup\{K_m:m\in [m_n,m_{n+1})\}$.
Since $\mathcal F$ is an ultrafilter, either $F_{even}=\bigcup\{[m_n,m_{n+1}): n\mbox{ is even}\}$
or $F_{odd}=\bigcup\{[m_n,m_{n+1}):n\mbox{ is odd}\}$ belongs to $\mathcal F$.
Without loss of generality, $F=F_{even}\in\mathcal F$.
For every $n\in\omega$ denote by $B_n$ the set $\bigcup_{m\in[m_{2n},m_{2n+1})}\bigcup_{k\in K_m}\cup v_k$
and note that $B_n$ is $u_n$-bounded. We claim that $\{B_n:n\in\omega\}$
is an $\omega$-cover of $X$. Indeed, given any finite subset $S$ of $X$,
for every $x\in S$ denote by $F_x$ the set $\{n\in\omega:I_s(x,w)\cap K_n\neq\emptyset\}$
and note that $F_x\in\mathcal F$. Since $\mathcal F$ is centered, there exists $m\in F\cap\bigcap_{x\in S}F_s$.
Let $n\in\omega$ be such that $m\in [m_{2n}, m_{2n+1})$. We claim that $S\subset B_n$.
Indeed, for every $x\in S$ there exists $k\in K_m\cap I_s(x,w)$, and thus $x\in\cup v_k\subset B_n$,
which finishes our proof.
 \end{proof}

Another application of Theorem~\ref{ch0} takes its origin in classical paper
\cite{Hu} of W.Hurewicz, where it was shown that a metrizable separable space
$X$ is Menger if and only if for every continuous function $f:X\rightarrow\mathbb R^\omega$ 
the image $f(X)$ is not dominating with respect to the eventual dominance preorder.
When $X$ is zero dimensonal, the same assertion holds for continuous functions
 $f:X\rightarrow\mathbb N^\omega$. 
Trying to generalize the above result 
 outside of metrizable separable spaces,
 all one can hope is the realm of Lindel\"of spaces
(every Menger topological space $X$ is obviously Lindel\"of). 
However, this obstacle may be overcome by restriction to countable covers
 in  the  definitions of the Menger  property. 
\begin{definition}
A topological space $X$ has the property $E^\ast_\omega$,
if for every sequence $(u_n)_{n\in\omega}$ of countable open covers of $X$
there exists a sequence $(v_n)_{n\in\omega}$ such that every $v_n$ is a finite subset
of $u_n$ and $\bigcup_{n\in\omega}\cup v_n=X$.  
\end{definition} 
It is clear that a topological space $X$ is Menger if and only if it has the property
$E^\ast_\omega$ and is Lindel\"of, and every countably compact noncompact space
has the property $E^\ast_\omega$ but fails to be Menger. 
The ideas of Hurewicz  still work for perfectly normal spaces: a perfectly normal
space $X$ has the property $E^\ast_\omega$ if and only if there is no continuous function
$f:X\rightarrow\mathbb R^\omega$ such that $f(X)$ is dominating,
see \cite{BH}. However, because of topological spaces $X$ such that
every continuous function $f:X\rightarrow\mathbb R$ is constant
the above characterization of the property $E^\ast_\omega$ is not true
for all topological spaces. For example, consider the topological space $Z=(\mathbb R,\tau)$,
where $\tau=\{(-\infty,a):a\in\mathbb R\}$. It is a simple exercise to show that
$Z^\omega$ is Lindel\"of and not Menger, but every continuous function $f:Z^\omega\rightarrow\mathbb R$
is constant. 

Theorem~\ref{ch0} enables us to prove a general characterization of the property $E^\ast_\omega$ 
involving compact-valued upper semicontinuous maps.
\begin{theorem}
A topological space $X$ has the property $E^\ast_\omega$ if and only if   $\Phi(X)\neq\mathbb N^\omega$ for every
compact-valued upper semicontinuous function $\Phi:X\Rightarrow\mathbb N^\omega$.
\end{theorem} 
\begin{proof}
 The ``only if'' part
follows from Lemma~\ref{mf}, which remains valid for the property $E^\ast_\omega$.
To prove the "if" part, we have to find an upper semicontinuous compact-valued surjective 
map $\Phi:X\Rightarrow\mathbb N^\omega$ provided $X$ has not  
the property $E^\ast_\omega$. A literal repetition of the proof of Lemma~\ref{crucial1}
gives us  a semifilter $\mathcal U\in\mathsf U(X)$
which fails to Menger. 
As it was shown in the proof of Corollary~\ref{ccrucial},
there exists a sequence $(\Phi_n)_{n\in\omega}$ of compact-valued upper semicontinuous maps 
from $X$ into $[\omega]^{\aleph_0}$ such that $\mathcal U=\bigcup_{n\in\omega}\Phi_n(X)$. Since the union of countably many
spaces with the Menger property is Menger, there exists $\Phi\in\{\Phi_n:n\in\omega\}$
such that the topological space $\Phi(X)$ has not the Menger property.   
Using the already mentioned result of Hurewicz, we can find a continuous map
$f:\Phi(X)\rightarrow\mathbb N^\omega$ such that $T=f(\Phi(X))$ is dominating in $\mathbb N^\omega$ with respect to $\leq^\ast$.
Next, we shall find a continuous map $g:T\rightarrow\mathbb N^\omega$ such that $g(T)$
is dominating in the following more strong sense:
for every $x\in\mathbb N^\omega$ there exists $y\in g(T)$ such that $y_n\geq x_n$ for all
$n\in\omega$. To find such the map $g$ it sufficies to note, that if none of the maps
$g_i:T\ni (x_n)_{n\in\omega}\mapsto(x_{n+i})_{n\in\omega}$ has this property, then
$T$ fails to be dominating.
And finally, consider the set-valued map $\Psi:g(T)\Rightarrow\mathbb N^\omega$,
$\Psi:g(T)\ni x\mapsto\{y\in\mathbb N^\omega:\forall n\in\omega (y_n\leq x_n)\}$. A direct verification shows that
$\Psi$ is  compact-valued and upper semicontinuous and $(\Psi\circ g\circ f\circ\Phi)(X)=\mathbb N^\omega$.
\end{proof}

\textbf{Acknowledgements.} The author wishes to express his thanks to prof.
Taras Banakh, who pointed out onto the important role of semifilters in selection
principles on topological spaces, and to prof. Boaz Tsaban, who  made several helpful
comments during the preparation of the paper.

Department of Mechanics and Mathematics, Ivan Franko Lviv National University,
Universytetska 1, Lviv, 79000, Ukraine.

\textit{E-mail address:} \texttt{lzdomsky@rambler.ru}


\end{document}